\documentclass[11pt,reqno]{amsart}
\usepackage{amsmath}
\usepackage{amssymb}
\usepackage{latexsym}

\newcommand{\Cplx}{\mathbb {C}}     
\newcommand{\halmos}{\rule{5pt}{5pt}}

\numberwithin{equation}{section}
\newtheorem{prop}{\bf Proposition}[section]
\newtheorem{thm}[prop]{\bf Theorem}
\newtheorem{df}[prop]{\bf Definition}

\newtheorem{cor}[prop]{\bf Corollary}

\newtheorem{rmk}[prop]{\bf Remark}

\begin{document}
\title[]
{Solutions to variant of $q$-hypergeometric equation of degree 2 associated with $q$-middle convolution}
\author[Yumi ARAI]{Yumi ARAI}
\address{Department of Mathematics, Ochanomizu University, 2-1-1 Otsuka, Bunkyo-ku, Tokyo 112-8610, Japan}
\email{araiyumi.math@gmail.com}
\subjclass[2020]{33D60, 33D15}
\keywords{$q$-hypergeometric equation, $q$-middle convolution, integral representation of solution}
\begin{abstract}
We investigate the integral representations of solutions to the variant of $q$-hypergeometric equation of degree 2  obtained through $q$-middle convolution by using transformation formulas for $q$-hypergeometric series.
We show the correspondence between these integral solutions and solutions obtained by other methods.
We also show the linear relationships among the integral solutions.
\end{abstract}
\maketitle

\section{Introduction}
In this paper, we describe the results of a study on the solutions to the variant of $q$-hypergeometric equation of degree 2 obtained by $q$-middle convolution.

The theory of middle convolution was introduced by Katz in his publication \textit{Rigid Local Systems} \cite{Katz}.
According to this theory, a Fuchsian differential equation is transformed into another Fuchsian differential equation by addition and middle convolution.
Subsequently, Dettweiler and Reiter reformulated Katz's theory in terms of linear algebra in \cite{DR2}. 
Here, an $n$-th order Fuchsian differential equation is a differential equation of the form
\begin{equation} \label{eq:FuchA}
\frac{d}{dx} Y(x) =\Bigl( \frac{A_1}{x-t_1}+\frac{A_2}{x-t_2}+\dots + \frac{A_r}{x-t_r} \Bigr) Y(x) ,
\end{equation}
where $Y(x)$ is a column vector with $n$ elements, and $A_1, A_2, \ldots, A_r$ are $n \times n$ constant matrices.
Set
\begin{equation*}
G_k = 
   \begin{pmatrix}     
       &        & O                      &        &  \\
   A_1 & \cdots & A_k + \lambda I_n & \cdots & A_r \\
       &        & O                      &        &  
    \end{pmatrix}
   (k), \quad 1\leq k \leq r, \quad \lambda \in \Cplx ,
\end{equation*}
where $I_n$ is the identity matrix of size $n$.
We call the correspondence of the tuples of matrices $(A_1, \ldots, A_r) \mapsto (G_1, \ldots, G_r)$ convolution.
The convolution is associated with Euler's integral transformation.
Let $Y(x)$ be a solution to equation (\ref{eq:FuchA}).
Set
\begin{equation*}
W_j(x) = \frac{Y(x)}{x - t_j},\ (j = 1, \ldots, r), \quad 
W(x) = \begin{pmatrix}
         W_1(x) \\
         \vdots \\
         W_r(x)
       \end{pmatrix} .
\end{equation*}
The function $Z(x) = \int_{\Delta} W(s) (x - s)^{\lambda} ds$ satisfies the Fuchsian differential equation
\begin{equation*}
\frac{d}{dx} Z(x) =\Bigl( \frac{G_1}{x-t_1}+\frac{G_2}{x-t_2}+\dots + \frac{G_r}{x-t_r} \Bigr) Z(x) ,
\end{equation*}
where $\Delta$ is an appropriate cycle (Euler's integral transformation).
We define two subspaces of $\Cplx^{Nm}$ as
\begin{equation*}
\mathcal{K} = \begin{pmatrix}
                \ker A_1 \\
                \vdots \\
                \ker A_r
              \end{pmatrix} , \quad
\mathcal{L} = \ker (G_1 + \cdots + G_r). 
\end{equation*}
Let $\overline{G}_k$ be the matrix induced by the action of $G_k$ on $\Cplx^{Nm}/(\mathcal{K} + \mathcal{L})$.
We call the correspondence of the tuples of matrices $(A_1, \ldots, A_r) \mapsto (\overline{G}_1, \ldots, \overline{G}_r)$ middle convolution.

Sakai and Yamaguchi constructed a $q$-analogue of the middle convolution in \cite{SY}.
The target of the $q$-middle convolution is a linear system of $q$-difference equations
\begin{equation}
E_{{\bold B}, {\bold b}}:Y(qx) = B(x)Y(x), \quad  B(x)=B_{\infty}+\sum_{i=1}^N \frac{B_i}{1-x/b_i},
\label{eq:EBb}
\end{equation}
where $Y(x)$ is a column vector with $m$ elements, $B_1, B_2, \ldots , B_N$ are $m \times m$ constant matrices, and $b_1, b_2, \ldots, b_N$ are non-zero complex numbers that are mutually distinct.

In this paper, we use the following notations:
\begin{align*}
&(a;q)_n= \prod_{j=0}^{n-1} (1-q^j a), \qquad (a_1, a_2, \dots, a_N;q)_n =(a_1;q)_n(a_2;q)_n \cdots (a_N;q)_n,\\
&(a;q)_{\infty}= \prod_{j=0}^{\infty} (1-q^j a), \qquad (a_1, a_2, \dots, a_N;q)_{\infty} =(a_1;q)_{\infty}(a_2;q)_{\infty} \cdots (a_N;q)_{\infty},\\
&{}_r \phi_{r-1} \left( \begin{array}{c} a_1, a_2, \ldots, a_r\\ b_1, \ldots, b_{r-1} \end{array} ;q, z \right)
= \sum_{n=0}^{\infty} \frac{( a_1, a_2, \ldots, a_r;q)_n}{(q;q)_n(b_1, \ldots, b_{r-1};q)_n} z^n,
\end{align*}
where $q$ is a complex number satisfying $0< |q| <1$.

Here, we recall the definitions of $q$-convolution and $q$-middle convolution and the theorem of $q$-integral transformation. 
\begin{df} [$q$-convolution] \cite[Definition 1.5]{SY} 
\label{def:qc} 
Let ${\bold B}= (  B_{\infty}; B_{1} ,\dots ,B_N ) $ be the tuple of $m\times m $ matrices and $\lambda \in \Cplx $.
Set $B_0 = I_m - B_\infty - B_{1} - \dots -B_N$.
We define the $q$-convolution $c_\lambda : (  B_{\infty};B_{1} ,\dots ,B_N ) \mapsto (  F_{\infty};F_{1} ,\dots ,F_N )$ as follows;
\begin{align}
&{\bold F} = ( F_\infty ; F_1, \dots , F_N) \mbox{ \rm  is a tuple of $(N+1)m \times (N+1)m$ matrices,}  \nonumber\\
&F_i = 
   \begin{pmatrix}     
       &        & O                      &        &  \\
   B_0 & \cdots & B_i - (1-q^\lambda)I_m & \cdots & B_N \\
       &        & O                      &        &  
    \end{pmatrix}
   (i+1), \quad 1\leq i \leq N,   \label{eq:bF}\\
&F_\infty = I_{(N+1)m} - \widehat{F}, \quad 
\widehat{F} =
    \begin{pmatrix}
        B_0 & \cdots & B_N \\
     \vdots & \ddots & \vdots \\
        B_0 & \cdots & B_N
    \end{pmatrix} \nonumber
\end{align}
\end{df}
The $q$-convolution in Definition \ref{def:qc} induces the correspondence of the linear $q$-difference equations
\begin{align*}
& E_{{\bold B}, {\bold b}}:Y(q x) = B(x)Y(x),\  B(x) = B_{\infty} + \sum^{N}_{i = 1}\frac{B_{i}}{1 - x/b_{i}} \\
&\qquad\qquad \mapsto \  E_{{\bold F}, {\bold b}}:\widehat{Y}(q x) = F(x)\widehat{Y}(x),\  F(x) = F_{\infty} + \sum^{N}_{i = 1}\frac{F_{i}}{1 - x/b_{i}} .
\end{align*}
Sakai and Yamaguchi introduced $q$-convolution replacing the Riemann integral in Dettweiler and Reiter's convolution with the infinite sum 
\begin{equation}
\label{qin}
\int_{0}^{\infty}f(s)d_q s := (1-q)\sum_{n=-\infty}^{\infty}q^nf(q^n),
\end{equation}
which is a $q$-analogue of Riemann integral.
In \cite{AT} with Takemura, equation (\ref{qin}) was extended as the following Jackson integral
\begin{equation} \label{qinxi}
 \int^{\xi \infty }_{0}f(s) \, d_{q}s = (1-q)\sum^{\infty}_{n=-\infty}q^{n} \xi f(q^n \xi ).
\end{equation}
Here, an additional parameter $\xi \in \mathbb{C} \setminus \{ 0 \}$ was introduced.
Furthermore, we investigated the convergence conditions of the integral representations of solutions.
In \cite{AT}, we extended Sakai and Yamaguchi's theorem related to $q$-integral transformation as follows:
\begin{thm} \cite[Theorem 2.6]{AT} \label{thm:qcintAT}
Let $Y(x)$ be a solution of $E_{{\bold B}, {\bold b}}$.
Then the function $\widehat{Y}(x) = {}^t (\widehat{Y}_{0}(x), \ldots, \widehat{Y}_{N}(x))$ defined by
\begin{align} \label{eq:qcintAT}
\begin{split}
& \widehat{Y}_{i}(x) = \int^{\xi \infty}_{0}\frac{P_{\lambda}(x, s)}{s-b_{i}}Y(s) \, d_{q}s, \; (i=0,\dots ,N), \\
& b_0 = 0, \quad P_{\lambda}(x,s) = \frac{(q^{\lambda+1}s/x;q)_{\infty}}{(qs/x;q)_{\infty}}
\end{split}
\end{align}
is convergent if the parameters satisfy certain conditions ((a) and (b) in \cite[Proposition 2.5]{AT}), and satisfies the equation $E_{{\bold F}, {\bold b}}$ for any $\xi \in \mathbb{C} \setminus \{ 0 \}$.
\end{thm}
The original theorem by Sakai and Yamaguchi is restricted to the case where $\xi = 1$ in the Jackson integral and does not consider convergence.
\begin{df}$(q$-middle convolution$)$ \cite[Definition 1.6]{SY} \label{def:qmc}
We define the $\mbox{\boldmath $F$}$-invariant subspaces $\mathcal{K}$ and $\mathcal{L} $ of $(\mathbb{C}^m )^{N+1}$ as follows;
\begin{equation} \label{eq:KL}
    \mathcal{K} = 
    \left( \hspace{-.5em}
        \begin{array}{c}
            \mathrm{ker}B_0 \\
            \vdots \\
            \mathrm{ker}B_N
        \end{array}
    \hspace{-.5em} \right),  \quad
    \mathcal{L} = 
    \mathrm{ker}(\widehat{F} - (1 - q^{\lambda})I_{(N+1)m}).
\end{equation}
We denote the action of $F_k$ on the quotient space $(\mathbb{C}^m )^{N+1}/(\mathcal{K} + \mathcal{L})$ by $\overline{F}_k$ $(k=\infty ,1, \dots ,N)$.
Then the $q$-middle convolution $mc_\lambda$ is defined by the correspondence $ E_{{\bold B}, {\bold b}}\mapsto E_{\overline{{\bold F}}, {\bold b}} $, where $ \overline{{\bold F}} = ( \overline{F}_{\infty}; \overline{F}_{1} ,\dots ,\overline{F}_N ) $.
\end{df}

In \cite{AT}, we applied $q$-convolution and $q$-middle convolution to specific $q$-difference equations and derived $q$-hypergeometric equation and its variants.
Additionally, we obtained integral solutions to these equations and their convergence conditions.
By applying the $q$-convolution to the $q$-difference equation satisfied by the function $y(x) = x^{\mu} (\alpha x; q)_{\infty}/(\beta x; q)_{\infty}$, we derived a single second-order $q$-difference equation. 
This equation corresponds to the standard form of $q$-hypergeometric equation 
\begin{equation} \label{eq:sta-qhg-s1}
(x-q)h(x/q)+(abx-c)h(qx)-\{(a+b)x-q-c \}h(x)=0 .
\end{equation}
The $q$-hypergoemetric series
\begin{equation*}
{}_2 \phi_1 \left( \begin{array}{c} a_1, a_2 \\ b \end{array} ;q, z \right)
= \sum_{n=0}^{\infty} \frac{( a_1, a_2 ;q)_n}{(q;q)_n(b ;q)_n} z^n  
\end{equation*}
introduced by Heine in 1846 (\cite{Hei}) satisfies this equation.
By taking the limit $q \to 1$, equation (\ref{eq:sta-qhg-s1}) tends to the hypergeometric differential equation
\begin{equation*}
 z(1 - z) \frac{d^2y}{dz^2} + (\gamma - (\alpha + \beta + 1)z) \frac{dy}{dz} - \alpha\beta y = 0. 
\end{equation*}
We also obtained the integral representations of solutions by Theorem \ref{thm:qcintAT}.
By specializing the value of the newly introduced parameter $\xi$, we could express the integral solutions in the form of $q$-hypergeometric series ${}_2 \phi_1$.

%
The function $y(x)=x^{\mu}(\alpha_1 x, \alpha_2 x;q)_{\infty}/(\beta_1 x, \beta_2 x ;q)_{\infty}$ satisfies the single $q$-difference equation $y(qx) = (B_{\infty} + B_1/(1 - \alpha_1 x) + B_2/(1 - \alpha_2 x)) y(x)$,  where $B_{\infty}$, $B_1$, $B_2$ are $1 \times 1$ matrices.
In \cite{AT}, we applied $q$-middle convolution to the tuple $(B_{\infty}; B_1, B_2)$ and derived the single second-order $q$-difference equations, which is correspond to the variant of $q$-hypergeometric equation of degree 2 
\begin{align}
& (x-q^{h_1 +1/2} t_1) (x - q^{h_2 +1/2} t_2) g(x/q) + q^{k_1 +k_2} (x - q^{l_1-1/2}t_1 ) (x - q^{l_2 -1/2} t_2) g(q x) \nonumber \\
&\hspace{80pt}  -[ (q^{k_1} +q^{k_2} ) x^2 +E x + p ( q^{1/2}+ q^{-1/2}) t_1 t_2 ] g(x) =0,  \label{eq:v2-intro}\\
&\hspace{60pt} p= q^{(h_1 +h_2 + l_1 + l_2 +k_1 +k_2 )/2 } , \quad E= -p \{ (q^{- h_2 }+q^{-l_2 })t_1 + (q^{- h_1 }+ q^{- l_1 }) t_2 \}. \nonumber
\end{align}
The variant of $q$-hypergeometric equation of degree 2 was introduced by Hatano, Matsunawa, Sato, and Takemura in \cite{HMST} in 2022.
This equation is a $q$-analogue of the second-order Fuchsian differential equation with three singularities $\{ t_1, t_2, \infty \}$, and is a special case of the $q$-Heun equation.
In the standard form of $q$-hypergeometric equation, the coefficients are linear polynomials in $x$, whereas in equation (\ref{eq:v2-intro}), they are quadratic polynomials.
The term ‘degree 2’ refers to the degree of the coefficients.
In \cite{AT}, we obtained the solutions to the variant of $q$-hypergeometric equation of degree 2 by the $q$-integral transformation associated with the $q$-middle convolution.
As in the case of the standard form of $q$-hypergeometric equation, if we substitute specific values to $\xi$, then the integral solutions are expressed by $q$-hypergeometric series $_3 \phi_2$.
A rigorous analysis of the convergence revealed that there are two cases depending on how the value of $\xi$ is chosen.
One case is when the integral solution satisfies the homogeneous variant of $q$-hypergeometric equation of degree 2, and the other case is when it satisfies the variant of $q$-hypergeometric equation of degree 2 with a non-homogeneous term (for details, see Section \ref{sec:is-vqhg2}). 

The solutions to equation (\ref{eq:v2-intro}) had already been obtained earlier in \cite{HMST, MST}, and additionally, by Fujii and Nobukawa in \cite{FN}.
In \cite{HMST, MST}, the solutions to equation (\ref{eq:v2-intro}) were obtained by substituting three formal solutions $g(x) = x^{-\rho} \sum_{k=0}^{\infty}c_{k} x^{-k}$, $g(x) = x^{\widetilde{\lambda}} \sum_{m=0}^{\infty}c_{m} (x/(q^{l_1-1/2} t_1) \\ ;q)_{m}$, and $g(x) = x^{-k_1} \sum_{n=0}^{\infty}c_{n} 
(q^{l_1-1/2} t_1/x;q)_{n}$ into equation (\ref{eq:v2-intro}) and determining the coefficients $c_{k}$, $c_{m}$, and $c_{n}$. 
In \cite{FN}, Fujii and Nobukawa first obtained the solutions to the variant of $q$-hypergeometric equation of degree 3 through three approaches.
First, they analyzed the equation in terms of its point configuration.
The second approach involved constructing solutions by integral representations.
Specifically, Jackson integrals of the Jordan-Pochhammer type are employed to derive explicit solutions.
By appropriately transforming these integrals and choosing different integral paths, they obtained several distinct solutions.
Furthermore, they transformed the integral solutions into $q$-hypergeometric series and obtained series solutions as well.
By taking the limit of some parameters in the series solutions to the variant of $q$-hypergeometric equation of degree 3, they obtained the solutions to the variant of $q$-hypergeometric equation of degree 2. 

We think that it is necessary to compare the integral solutions obtained by the $q$-middle convolution with the already found solutions to see if there is a correspondence between them.

Since the variants of $q$-hypergeometric equation were only recently introduced, there is still little known about its solutions. 
In this study, we aim to elucidate the properties of the solutions, focusing in particular on investing the linear relationships among them.
For the purpose of our research, we use several transformation formulas written in \textit{Basic Hypergeometric Series} (G.Gasper, M.Rahman) \cite{GR} as tools.

If we appropriately take the limit as $q \to 1$ in $q$-difference equations and their solutions, they transition to differential equations and their solutions.
By taking the limit of the results obtained by the $q$-middle convolution, it is also expected that these results can be applied to the study of special functions.
%

This paper is organized as follows. 
In Section \ref{sec:is-vqhg2}, as preliminary, we review the derivation of the variant of $q$-hypergeometric equation of degree 2 and its integral representation of solutions by applying $q$-middle convolution. 
In Section \ref{sec:correso}, we compare the integral solutions obtained through the $q$-middle convolution with the solutions found by different methods and give their correspondence.
In section \ref{sec:LinRel}, we discuss the linear relationships among the integral solutions obtained by the $q$-middle convolution.

\section{Integral representation of solutions to variant of $q$-hypergeometric equation of degree 2} \label{sec:is-vqhg2}

In \cite{AT}, Takemura and the author applied the $q$-middle convolution to the linear $q$-difference equation satisfied by the function 
\begin{equation*}
y(x)=(\alpha_1 x, \alpha_2 x;q)_{\infty}/(\beta_1 x, \beta_2 x;q)_{\infty} ,
\end{equation*}
derived the single second-order $q$-difference equation corresponding to the variant of $q$-hypergeometric equation of degree 2, and obtained its integral representations of solutions.
%
%
We review the content of \cite[Section 4.1.1]{AT}.

The function $y(x) = (\alpha_1 x, \alpha_2 x;q)_{\infty}/(\beta_1 x, \beta_2 x;q)_{\infty}$ satisfies the linear $q$-difference equation $E_{{\bold B}, {\bold b}}: y(qx) = B(x)y(x)$, where
\begin{equation}
\begin{split}
&B(x) = B_{\infty} + \frac{B_1}{1 - x/b_1} + \frac{B_2}{1 - x/b_2}, \ b_1 = \frac{1}{\alpha_1}, \ b_2 = \frac{1}{\alpha_2}, \\
&B_{\infty} = \frac{\beta_1\beta_2}{\alpha_1\alpha_2}, \ B_1 = \frac{(\alpha_1 - \beta_1)(\alpha_1 - \beta_2)}{\alpha_1(\alpha_1 - \alpha_2)}, \ B_2 = \frac{(\alpha_2 - \beta_1)(\alpha_2 - \beta_2)}{\alpha_2(\alpha_2 - \alpha_1)}.  
\end{split} \label{eq:Bb1b2}
\end{equation}
Note that $B_0 = 0$.
We apply the $q$-convolution to the tuple of matrices $(B_{\infty}; B_1, B_2)$ and set $c_{\lambda}(B_{\infty}; B_1, B_2) = (F_{\infty}; F_1, F_2)$.
Then we obtain the matrices
\begin{align*}
&F_1 = \begin{pmatrix}
        0   & 0                     & 0 \\
        0 & B_1 - 1 + q^{\lambda} & B_2 \\
        0   & 0                     & 0 
\end{pmatrix}, \ 
F_2 = \begin{pmatrix}
        0   & 0   & 0 \\
        0   & 0   & 0 \\       
        0 & B_1 & B_2 - 1 + q^{\lambda}
\end{pmatrix}, \\
&F_{\infty} = \begin{pmatrix}
               1 & -B_1    & -B_2 \\
               0 & 1 - B_1 & -B_2 \\       
               0 & -B_1    & 1 - B_2
\end{pmatrix}.
\end{align*}
The corresponding $q$-difference equation $E_{{\bold F}, {\bold b}}$ is written as
\begin{equation} \label{eq:K1F2}
\widehat{Y}(qx) = \left( F_{\infty} + \frac{F_1}{1 - \alpha_1 x} + \frac{F_2}{1 - \alpha_2 x} \right) \widehat{Y}(x) , \quad  
\widehat{Y}(x) = \begin{pmatrix}
                   \widehat{y}_0(x) \\
                   \widehat{y}_1(x) \\
                   \widehat{y}_2(x) 
                 \end{pmatrix}.
\end{equation}
In this case, the subspace $\mathcal{K}$ in Definition \ref{def:qmc} is spanned by the vector ${}^t(1,0,0)$.
On the other hand, $\mathcal{L} = \{ 0 \}$.
Hence, we have $\dim \mathcal{K} = 1$ and $\dim \mathcal{L} = 0$.
From equation (\ref{eq:K1F2}), we can obtain 
\begin{align} \label{eq:haty1y2}
\begin{split}  
&\begin{pmatrix} 
  \widehat{y}_1(qx) \\
  \widehat{y}_2(qx) 
\end{pmatrix}
= \left( \overline{F}_{\infty} + \frac{\overline{F}_1}{1 - \alpha_1 x} + \frac{\overline{F}_2}{1 - \alpha_2 x} \right)
\begin{pmatrix} 
  \widehat{y}_1(x) \\
  \widehat{y}_2(x) 
\end{pmatrix}, \\
&\overline{F}_1 = \begin{pmatrix}
                     B_1 - 1 + q^{\lambda} & B_2 \\
                     0                     & 0 
                   \end{pmatrix} , \ 
\overline{F}_2 =  \begin{pmatrix}
                     0   & 0 \\       
                     B_1 & B_2 - 1 + q^{\lambda}
                   \end{pmatrix} , \ 
\overline{F}_{\infty} = \begin{pmatrix}
                          1 - B_1 & -B_2 \\       
                          -B_1    & 1 - B_2
                        \end{pmatrix}.
\end{split}
\end{align}
The matrices $\overline{F}_1$, $\overline{F}_2$, and $\overline{F}_{\infty}$ are the lower-right $2 \times 2$ submatrices of $F_1$, $F_2$, and $F_{\infty}$, and they can be regarded as those obtained by the $q$-middle convolution.
By eliminating $\widehat{y}_2(x)$ from eq. (\ref{eq:haty1y2}), we obtain the single second-order $q$-difference equation for $\widehat{y}_1(x)$
\begin{align}
&\Bigl( x-\frac{q^{\lambda+1}}{\beta_1} \Bigr) \Bigl( x-\frac{q^{\lambda+1}}{\beta_2} \Bigr) \widehat{y}_1(x/q) + \frac{\alpha_1\alpha_2}{\beta_1\beta_2} \Bigl( x-\frac{1}{\alpha_1} \Bigr) \Bigl(x-\frac{q}{\alpha_2} \Bigr) \widehat{y}_1(qx) \notag \\
 &- \Bigl[ \Bigl( \frac{\alpha_1\alpha_2}{\beta_1\beta_2}+1 \Bigr) x^2-\Bigl\{ q \Bigl( \frac{1}{\beta_1}+\frac{1}{\beta_2}\Bigr) +q^{\lambda}\frac{q\alpha_1+\alpha_2}{\beta_1\beta_2} \Bigr\} x+\frac{q^{\lambda+1}(1+q)}{\beta_1\beta_2} \Bigr] \widehat{y}_1(x) = 0. \label{eq:g1eqdeg2k}
\end{align}
Equation (\ref{eq:g1eqdeg2k}) is a special case of the variant of $q$-hypergeometric equation of degree 2 
\begin{equation}
\begin{split}
& (x-q^{h_1 +1/2} t_1) (x - q^{h_2 +1/2} t_2) g(x/q) + q^{k_1 +k_2} (x - q^{l_1-1/2}t_1 ) (x - q^{l_2 -1/2} t_2) g(q x) \\
&  -[ (q^{k_1} +q^{k_2} ) x^2 +E x + p ( q^{1/2}+ q^{-1/2}) t_1 t_2 ] g(x) =0, \label{eq:varqhgdeg2}\\
& p= q^{(h_1 +h_2 + l_1 + l_2 +k_1 +k_2 )/2 } , \quad E= -p \{ (q^{- h_2 }+q^{-l_2 })t_1 + (q^{- h_1 }+ q^{- l_1 }) t_2 \} 
\end{split}
\end{equation}
with the constraint $k_2=0$.
\begin{prop}$($\cite[Proposition 4.1]{AT}$)$  \label{prop:paramrel}
Assume that $\widehat{y}_1(x)$ satisfies equation (\ref{eq:g1eqdeg2k}).
Set $g(x) = x^{-k_2} \widehat{y}_1(x)$ and
\begin{align*}
& \alpha_1 = \frac{q^{- l_1 + 1/2}}{t_1}, \quad \alpha_2  = \frac{q^{-l_2 +3/2} }{t_2 } ,\quad q^{\lambda }=q^{(h_1+h_2-l_1-l_2-k_1 +k_2 +1)/2} \\
& \beta_1  = \frac{q^{ (-h_1+h_2-l_1-l_2-k_1 +k_2)/2 +1} }{ t_1} \left(= \frac{q^{\lambda -h_1 +1/2}}{t_1}\right),\\
& \beta_2 = \frac{q^{(h_1- h_2-l_1-l_2-k_1 +k_2)/2 +1}}{ t_2 } \left(= \frac{q^{\lambda - h_2 + 1/2}}{t_2}\right).  
\end{align*}
Then $g(x)$ satisfies the variant of $q$-hypergeometric equation of degree 2 given by (\ref{eq:varqhgdeg2}).
\end{prop}

Next, we consider the integral representation of solutions.
In this case, the assumption regarding the convergence conditions in Theorem \ref{thm:qcintAT} does not hold, and thus Theorem \ref{thm:qcintAT} cannot be applied.
However, we represent the formal solution as
\begin{align}
\widehat{y}_1(x) 
&= \int_0^{\xi \infty} \frac{P_{\lambda}(x,s)}{s-1/\alpha_1} y(s) d_q s  \nonumber \\
&= (q-1)\alpha_1 \sum_{n=-\infty}^{\infty} q^n \xi \frac{(q^{\lambda+n+1}\xi/x, q^{n+1}\xi\alpha_1, q^n\xi\alpha_2;q)_{\infty}}{(q^{n+1}\xi/x, q^n\xi\beta_1, q^n\xi\beta_2;q)_{\infty}}. 
\label{eq:yhata1}
\end{align}
By substituting $\xi=1/\alpha_1$, $\xi=1/\alpha_2$, or $\xi=q^{-\lambda}x$ into equation (\ref{eq:yhata1}), it is written as ${}_3\phi_2$ series.
Let $y_{\alpha_1}(x)$, $y_{\alpha_2}(x)$, and $y_{\lambda}(x)$ denote the functions obtained by substituting $\xi=1/\alpha_1$, $\xi=1/\alpha_2$, and $\xi=q^{-\lambda}x$, respectively:
\begin{align}
& y_{\alpha_1}(x) = (q-1)\frac{(q^{\lambda+1}/(\alpha_1 x),q,\alpha_2/\alpha_1;q)_{\infty}}{(q/(\alpha_1 x),\beta_1/\alpha_1,\beta_2/\alpha_1;q)_{\infty}} \notag \\
&\hspace{150pt}\cdot\!_3\phi_2 \left( \begin{array}{c} q/(\alpha_1 x), \beta_1/\alpha_1,\beta_2/\alpha_1 \\ q^{\lambda+1}/(\alpha_1 x), \alpha_2/\alpha_1  \end{array} ;q, q \right), \label{eq:ya1m0} \\
& y_{\alpha_2}(x) = (q-1)q\frac{\alpha_1}{\alpha_2}\frac{(q^{\lambda+2}/(\alpha_2 x),q^2\alpha_1/\alpha_2,q;q)_{\infty}}{(q^2/(\alpha_2 x),q\beta_1/\alpha_2,q\beta_2/\alpha_2;q)_{\infty}} \notag\\
&\hspace{150pt}\cdot\!_3\phi_2 \left( \begin{array}{c} q^2/(\alpha_2 x), q\beta_1/\alpha_2,q\beta_2/\alpha_2 \\ q^{\lambda+2}/(\alpha_2 x),q^2\alpha_1/\alpha_2  \end{array} ;q, q \right), \label{eq:ya2m0} \\
& y_{\lambda}(x) = (q-1)q^{-\lambda}\alpha_1 x \frac{(q,q^{-\lambda+1}\alpha_1 x,q^{-\lambda}\alpha_2 x;q)_{\infty}}{(q^{-\lambda+1},q^{-\lambda}\beta_1 x,q^{-\lambda}\beta_2 x;q)_{\infty}} \notag\\
&\hspace{150pt}\cdot\!_3\phi_2 \left( \begin{array}{c} q^{-\lambda+1},q^{-\lambda}\beta_1 x,q^{-\lambda}\beta_2 x \\ q^{-\lambda+1}\alpha_1 x,q^{-\lambda}\alpha_2 x  \end{array} ;q, q \right). \label{eq:eyalm0}
\end{align}
In \cite[Proposition 4.3]{AT}, it was shown that functions (\ref{eq:ya1m0})-(\ref{eq:eyalm0}) converge and satisfy the non-homogeneous version of equation (\ref{eq:g1eqdeg2k}) given by
\begin{align}
&\Bigl( x-\frac{q^{\lambda+1}}{\beta_1} \Bigr) \Bigl( x-\frac{q^{\lambda+1}}{\beta_2} \Bigr) \widehat{y}_1(x/q) + \frac{\alpha_1\alpha_2}{\beta_1\beta_2} \Bigl( x-\frac{1}{\alpha_1} \Bigr) \Bigl( x-\frac{q}{\alpha_2} \Bigr) \widehat{y}_1(qx) \nonumber \\
  &- \Bigl[ \Bigl( \frac{\alpha_1\alpha_2}{\beta_1\beta_2}+1 \Bigr) x^2- \Bigl\{ q \Bigl( \frac{1}{\beta_1}+\frac{1}{\beta_2} \Bigr) +q^{\lambda}\frac{q\alpha_1+\alpha_2}{\beta_1\beta_2} \Bigr\} x+\frac{q^{\lambda+1}(1+q)}{\beta_1\beta_2} \Bigr] \widehat{y}_1(x) \nonumber \\
  &-q(1-q)(1-q^{\lambda})\frac{\alpha_1}{\beta_1\beta_2}x  = 0.
\label{eq:y1nonhom}
\end{align}
For functions (\ref{eq:ya1m0})-(\ref{eq:eyalm0}), the difference between any two of these functions satisfies the homogeneous equation given by (\ref{eq:g1eqdeg2k}).
Let us also consider other solutions.
To substitute $\xi=1/\beta_1$, $\xi=1/\beta_2$, or $\xi=x$ into $\widehat{y}_1(x)$, we replace the functions $P_{\lambda}(x,s)$ and $y(s)$ in equation (\ref{eq:yhata1}) with
\begin{equation}
\widetilde{P}_{\lambda}(x,s)=(x/s)^{\lambda}\frac{(x/s;q)_{\infty}}{(q^{-\lambda}x/s;q)_{\infty}}, \qquad \widetilde{y}(s) = s^{\mu'}\frac{(q/(\beta_1 s), q/(\beta_2 s);q)_{\infty}}{(q/(\alpha_1 s), q/(\alpha_2 s);q)_{\infty}} \label{eq:Ply}
\end{equation}
with the condition $q^{\mu'}\alpha_1\alpha_2/\beta_1\beta_2 = 1$.
The function $P_{\lambda}(x,s)$ satisfies the $q$-difference equation $P_{\lambda}(qx,s) = (x - q^{\lambda}s) P_{\lambda}(x,s)/(x - s)$, and $\widetilde{P}_{\lambda}(x,s)$ also satisfies this equation.
The function $\widetilde{y}(s)$ satisfies the $q$-difference equation $\widetilde{y}(qx) = B(x)\widetilde{y}(x)$, where $B(x)$ is given by equation (\ref{eq:Bb1b2}).
The function $\widehat{y}_1(x)$ is rewritten as
\begin{align*}
\widehat{y}_1(x) = (1-q)x^{\lambda}\sum_{n=-\infty}^{\infty} (q^n \xi)^{\mu'-\lambda} \frac{(xq^{-n}/\xi, q^{1-n}/(\beta_1\xi), q^{1-n}/(\beta_2\xi); q)_{\infty}}{(xq^{-\lambda-n}/\xi, q^{-n}/(\alpha_1\xi), q^{1-n}/(\alpha_2\xi); q)_{\infty}}.
\end{align*}
Let $y_{\beta_1}(x)$, $y_{\beta_2}(x)$, and $y_x(x)$ denote the functions obtained by substituting $\xi=1/\beta_1$, $\xi=1/\beta_2$, and $\xi=x$ into $\widehat{y}_1(x)$, respectively:
\begin{align}
& y_{\beta_1}(x) = (1-q)\beta_1^{\lambda-\mu'}x^{\lambda}\frac{(\beta_1 x, q, q\beta_1/\beta_2 ;q)_{\infty}}{(q^{-\lambda}\beta_1 x, \beta_1/\alpha_1, q\beta_1/\alpha_2 ;q)_{\infty}} \notag\\
&\hspace{150pt}\cdot\!_3\phi_2 \left( \begin{array}{c} q^{-\lambda}\beta_1 x, \beta_1/\alpha_1,q\beta_1/\alpha_2 \\ \beta_1 x, q\beta_1/\beta_2  \end{array} ;q, q^{\lambda}\frac{\alpha_1\alpha_2}{\beta_1\beta_2} \right), \label{eq:yb1m0} \\
& y_{\beta_2}(x) = (1-q)\beta_2^{\lambda-\mu'}x^{\lambda}\frac{(\beta_2 x, q, q\beta_2/\beta_1 ;q)_{\infty}}{(q^{-\lambda}\beta_2 x,\beta_2/\alpha_1, q\beta_2/\alpha_2 ;q)_{\infty}} \notag\\
&\hspace{150pt}\cdot\!_3\phi_2 \left( \begin{array}{c} q^{-\lambda}\beta_2 x, \beta_2/\alpha_1,q\beta_2/\alpha_2 \\ \beta_2 x, q\beta_2/\beta_1  \end{array} ;q, q^{\lambda}\frac{\alpha_1\alpha_2}{\beta_1\beta_2} \right), \label{eq:yb2m0} \\
& y_x(x) = (1-q)q^{\lambda-\mu'}x^{\mu'} \frac{(q^2/(\beta_1 x), q^2/(\beta_2 x), q ;q)_{\infty}}{(q/(\alpha_1 x), q^2/(\alpha_2 x), q^{-\lambda+1} ;q)_{\infty}} \notag\\
&\hspace{150pt}\cdot\!_3\phi_2 \left( \begin{array}{c} q^{-\lambda+1}, q/(\alpha_1 x),q^2/(\alpha_2 x) \\ q^2/(\beta_1 x), q^2/(\beta_2 x)  \end{array} ;q, q^{\lambda}\frac{\alpha_1\alpha_2}{\beta_1\beta_2} \right). \label{eq:yxm0}
\end{align}
In \cite[Proposition 4.4]{AT}, it was shown that functions (\ref{eq:yb1m0})-(\ref{eq:yxm0}) converge under the condition $|q^{\lambda}\alpha_1\alpha_2/(\beta_1\beta_2)|<1$ and satisfy equation (\ref{eq:g1eqdeg2k}).

\section{Correspondence with already found solutions} \label{sec:correso}
\subsection{Correspondence with $q$-Appell series solution found in \cite{HMST, MST}}$ $ \label{subsec:qAp} 

In \cite{HMST} and \cite{MST}, three solutions to equation (\ref{eq:varqhgdeg2}) are given, one of which is expressed using the $q$-Appell series.
The $q$-Appell series is defined as
\begin{equation}
\Phi^{(1)}(a; b,b'; c;q; y, z) 
= \sum_{m=0}^{\infty}\sum_{n=0}^{\infty} \frac{(a ;q)_{m+n} (b; q)_m (b'; q)_n}{(c ;q)_{m+n} (q; q)_m (q; q)_n} y^m z^n.
\end{equation}
It is the $q$-analogue of the well-known Appell series
\begin{equation*}
F_1(a; b,b'; c;q; y, z) = \sum_{m=0}^{\infty} \sum_{n=0}^{\infty} \frac{(a)_{m+n} (b)_m (b')_n}{(c)_{m+n} m! n!} y^m z^n, \quad (a)_n = a(a+1) \cdots (a+n-1)
\end{equation*}
which is a two-variable extension of Gauss' hypergeometric series ${}_2F_1$.

To the $q$-hypergeometric series $_3\phi_2 (q^{-\lambda+1}, q/(\alpha_1 x), q^2/(\alpha_2 x); q^2/(\beta_1 x), q^2/(\beta_2 x);q,
\\q^{\lambda}\alpha_1\alpha_2/(\beta_1\beta_2))$ in function (\ref{eq:yxm0}), we apply Andrews' formula (\cite{And}, \cite[formula (10.3.4)]{GR})
\begin{equation}
\Phi^{(1)}(q^a; q^b,q^{b^{\prime}}; q^c;q; y,z)=\frac{(q^a, yq^b, zq^{b^{\prime}}; q)_{\infty}}{(q^c,y,z; q)_{\infty}}{}_3\phi_2 \left( \begin{array}{c} q^{c-a},y,z\\ yq^b,zq^{b^{\prime}} \end{array} ;q,q^a \right). \label{qappell}
\end{equation} 
Set
\begin{align*}
&q^{c-a}=q^{-\lambda+1},\quad y=q/(\alpha_1 x),\quad z=q^2/(\alpha_2 x), \notag \\
&yq^b=q^2/(\beta_1 x),\quad zq^{b^{\prime}}=q^2/(\beta_2 x),\notag \\
&q^a=q^{\lambda}\alpha_1\alpha_2/(\beta_1\beta_2).
\end{align*}
From $q^b=q\alpha_1/\beta_1$, $q^{b^{\prime}}=\alpha_2/\beta_2$, and $q^c=q\alpha_1\alpha_2/(\beta_1\beta_2)$, we obtain
\begin{align*}
&_3 \phi_2 \left( \begin{array}{c} q^{-\lambda+1}, q/(\alpha_1 x), q^2/(\alpha_2 x)\\q^2/(\beta_1 x), q^2/(\beta_2 x)  \end{array} ;q, q^{\lambda}\frac{\alpha_1\alpha_2}{\beta_1\beta_2}\right)\\
&= \frac{(q\alpha_1\alpha_2/(\beta_1\beta_2), q/(\alpha_1 x), q^2/(\alpha_2 x); q)_{\infty}}{(q^{\lambda}\alpha_1\alpha_2/(\beta_1\beta_2), q^2/(\beta_1 x), q^2/(\beta_2 x); q)_{\infty}}\\ 
&\qquad \cdot \Phi^{(1)}(q^{\lambda}\alpha_1\alpha_2/(\beta_1\beta_2); q\alpha_1/\beta_1, \alpha_2/\beta_2; q\alpha_1\alpha_2/(\beta_1\beta_2); q; q/(\alpha_1 x), q^2/(\alpha_2 x))\\
&= \frac{(q\alpha_1\alpha_2/(\beta_1\beta_2), q/(\alpha_1 x), q^2/(\alpha_2 x); q)_{\infty}}{(q^{\lambda}\alpha_1\alpha_2/(\beta_1\beta_2), q^2/(\beta_1 x), q^2/(\beta_2 x); q)_{\infty}}\\
&\qquad \cdot \sum_{n=0}^{\infty}\sum_{m=0}^{\infty} \frac{(q^{\lambda}\alpha_1\alpha_2/(\beta_1\beta_2) ;q)_{m+n} (q\alpha_1/\beta_1; q)_m (\alpha_2/\beta_2; q)_n}{(q\alpha_1\alpha_2/(\beta_1\beta_2) ;q)_{m+n} (q; q)_m (q; q)_n}
(q/(\alpha_1 x))^m (q^2/(\alpha_2 x))^n,
\end{align*}
which leads to
\begin{align}
y_x(x) &= (1-q)q^{\lambda-\mu'}\frac{(q, q\alpha_1\alpha_2/(\beta_1\beta_2); q)_{\infty}}{(q^{-\lambda+1}, q^{\lambda}\alpha_1\alpha_2/(\beta_1\beta_2); q)_{\infty}} \notag \\
& \cdot x^{\mu'}\sum_{n=0}^{\infty}\sum_{m=0}^{\infty} \frac{(q^{\lambda}\alpha_1\alpha_2/(\beta_1\beta_2) ;q)_{m+n} (q\alpha_1/\beta_1; q)_m (\alpha_2/\beta_2; q)_n}{(q\alpha_1\alpha_2/(\beta_1\beta_2) ;q)_{m+n} (q; q)_m (q; q)_n}
(q/(\alpha_1 x))^m (q^2/(\alpha_2 x))^n. \label{eq:yxqappell}
\end{align}
From Proposition \ref{prop:paramrel} and equation (\ref{eq:yxqappell}), we obtain the following proposition.
\begin{prop}
The function $y_x(x)$ corresponds to a constant multiple of the $q$-Appell series solution to equation (\ref{eq:varqhgdeg2}) given by
\begin{align} 
g(x) &= x^{-k_1}\sum_{n=0}^{\infty}\sum_{m=0}^{\infty} \frac{(q^{\widetilde{\lambda} +k_1}; q)_{m+n} (q^{\widetilde{\lambda} +k_1-h_2+l_2}; q)_m (q^{\widetilde{\lambda} +k_1-h_1+l_1}; q)_n}{(q^{k_1-k_2+1}; q)_{m+n} (q; q)_m (q; q)_n} \notag \\
&\hspace{150pt} \cdot (q^{l_1+1/2}t_1 x^{-1})^m (q^{l_2+1/2}t_2 x^{-1})^n, \label{eq:g1appell}
\end{align}
where $\mu'=-k_1+k_2$ and $\lambda=\widetilde{\lambda}+k_2$.
\end{prop}

\subsection{Correspondence with solutions obtained by Fujii and Nobukawa}$ $ \label{subsec:solFN}

Fujii and Nobukawa also investigated the solutions to the homogeneous variant of $q$-hypergeometric equation of degree 2 (\cite{FN}). 
Now we investigate the correspondence between our results and theirs.

In \cite{FN}, the variant of $q$-hypergeometric equation of degree 2 is written as $\mathcal{E}_2 f(x)=0$, where
\begin{align*}
\mathcal{E}_2 = [x^2(1-q^{\alpha}T_x)(B-AT_x) &- x(e_1(a)-q^{\alpha}e_1(b)T_x)(1-T_x) \notag \\ 
&\hspace{50pt} + e_2(a)B^{-1}(1-q^{-1}T_x)(1-T_x)]T_x^{-1}, \\
& Aa_1a_2 = q^{\alpha + 1}Bb_1b_2
\end{align*}
(\cite[equation (3.43)]{FN}). 
The notation $e_i$ denotes the elementary symmetric polynomial of degree $i$.
The equation $\mathcal{E}_2 f(x)=0$ is rewritten as
\begin{align}
&\biggl(x-\frac{a_1}{B}\biggr) \biggl(x-\frac{a_2}{B} \biggr) f(x/q) + q^{\alpha}\frac{A}{B} \biggl(x-\frac{b_1}{A}\biggr)\biggl(x-\frac{b_2}{A}\biggr) f(qx) \notag \\
&\qquad - \biggl[ \biggl(\frac{A}{B}+q^{\alpha}\biggr)x^2 - \biggl(\frac{a_1+a_2}{B}+q^{\alpha}\frac{b_1+b_2}{B}\biggr)x + \frac{a_1a_2}{qB^2}(1+q) \biggr] f(x)=0 \label{eq:e2}
\end{align}
and is transformed into equation (\ref{eq:varqhgdeg2}) with
\begin{align}
g(x) = x^{\lambda_0} f(x), \quad \alpha = \lambda_0 + k_1, \quad B/A = q^{-k_2 - \lambda_0}, \quad a_i/B = q^{h_i+1/2}t_i, \quad b_i/A = q^{l_i-1/2}t_i , \label{param:fg}
\end{align}
where $\lambda_0 = (h_1 + h_2 - l_1 - l_2 - k_1 - k_2 + 1)/2$
(\cite[equation (3.44)]{FN}). 
Assume that $\widehat{y}_1(x)$, $g(x)$, and $f(x)$ satisfy equations (\ref{eq:g1eqdeg2k}), (\ref{eq:varqhgdeg2}), and (\ref{eq:e2}), respectively.
Then, by Proposition \ref{prop:paramrel} and equation (\ref{param:fg}) the following relations are obtained:
\begin{equation}
\begin{split}
&f(x) = x^{-\lambda}\widehat{y}_1(x),\\
&\frac{a_1}{B} = \frac{q^{\lambda+1}}{\beta_1}, \quad \frac{a_2}{B} = \frac{q^{\lambda+1}}{\beta_2},\quad \frac{b_1}{A} = \frac{1}{\alpha_1}, \quad \frac{b_2}{A} = \frac{q}{\alpha_2}, \\
&\frac{A}{B} = q^{\lambda_0 + k_2} = q^{\lambda}, \ q^{\alpha} = q^{\lambda_0 + k_1} = q^{\lambda+k_1-k_2} = q^{\lambda}\frac{\alpha_1\alpha_2}{\beta_1\beta_2}. 
\label{param:fyhat}
\end{split}
\end{equation}
In \cite{FN}, Fujii and Nobukawa introduced the following functions for solutions to $\mathcal{E}_2 f(x)=0$ :
\begin{align}
&{}_3\phi_2 \left( \begin{array}{c} q^{\alpha}, A/B, a_1/(Bx) \\ A a_1/(B b_1), A a_1/(B b_2) \end{array} ;q, qB x/a_2 \right), \label{eq:FN432} \\
&\frac{(qAx/a_2;q)_{\infty}}{(qBx/a_2;q)_{\infty}} {}_3\phi_2 \left( \begin{array}{c} qb_1/a_2, A/B, Ax/b_1 \\ A a_1/(B b_1), qAx/a_2 \end{array} ;q, a_1/b_2 \right), \label{eq:FN433} \\
&\frac{(qAx/a_2;q)_{\infty}}{(qBx/a_2;q)_{\infty}} {}_3\phi_2 \left( \begin{array}{c} qb_1/a_2, qb_2/a_2, A/B \\ q^{1-\alpha}A/B , qAx/a_2 \end{array} ;q, q \right), \label{eq:FN434} \\
&\frac{(Ax/b_2;q)_{\infty}}{(b_1Bx/(a_1a_2);q)_{\infty}} {}_3\phi_2 \left( \begin{array}{c} q^{\alpha}, a_1/b_1, a_2/b_1 \\ a_1a_2/(b_1Bx) , a_1a_2/(b_1b_2) \end{array} ;q, qb_1/(Ax) \right), \label{eq:FN435} \\
&\frac{(Ax/a_2, qBx/a_1;q)_{\infty}}{(Ax/b_1, Ax/b_2;q)_{\infty}} {}_3\phi_2 \left( \begin{array}{c} qB/A, a_2/b_2, qb_1/(Ax) \\ q^{1-\alpha}a_2/b_2 , qa_2/(Ax) \end{array} ;q, q \right), \label{eq:FN436} \\
&{}_3\phi_2 \left( \begin{array}{c} q^{\alpha}, A/B, Ax/b_1 \\ A a_1/(B b_1), A a_2/(B b_1) \end{array} ;q, q \right), \label{eq:FN437}
\end{align}
where $a_1a_2A = q^{\alpha+1}b_1b_2 B$ (\cite[Theorem 4.2]{FN}).
After modifying functions (\ref{eq:FN433}), (\ref{eq:FN435}), and (\ref{eq:FN436}) as
\begin{align*}
&(\ref{eq:FN433}) : {}_3\phi_2 \left( \begin{array}{c} qb_1/a_2, A/B, Ax/b_1 \\ A a_1/(B b_1), qAx/a_2 \end{array} ;q, a_1/b_2 \right)
\rightarrow
{}_3\phi_2 \left( \begin{array}{c} qb_2/a_2, A/B, Ax/b_1 \\Aa_1/(Bb_1), qAx/a_2 \end{array} ;q, a_1/b_2 \right),\\
&(\ref{eq:FN435}) : \displaystyle \frac{(Ax/b_2;q)_{\infty}}{(b_1Bx/(a_1a_2);q)_{\infty}} 
\rightarrow 
\displaystyle \frac{(Ax/b_2;q)_{\infty}}{(qb_1Bx/(a_1a_2);q)_{\infty}},\\
&(\ref{eq:FN436}) : \displaystyle \frac{(Ax/a_2, qBx/a_1;q)_{\infty}}{(Ax/b_1, Ax/b_2;q)_{\infty}} 
\rightarrow 
\displaystyle \frac{(Ax/b_1, Ax/b_2;q)_{\infty}}{(Ax/a_2, qBx/a_1;q)_{\infty}},
\end{align*}
we rewrite functions (\ref{eq:FN432})-(\ref{eq:FN437}) by relations (\ref{param:fyhat}) and refer to them as $f_{32}(x)$-$f_{37}(x)$:
\begin{align}
&f_{32}(x) = x^{\lambda}{}_3\phi_2 \left( \begin{array}{c} q^{\lambda}\alpha_1\alpha_2/(\beta_1\beta_2), q^{\lambda}, q^{\lambda+1}/(\beta_1 x) \\ q^{\lambda+1}\alpha_1/\beta_1, q^{\lambda}\alpha_2/\beta_1 \end{array} ;q, q^{-\lambda}\beta_2 x \right) \label{eq:fn432-1},\\
&f_{33}(x) = x^{\lambda}\frac{(\beta_2 x;q)_{\infty}}{(q^{-\lambda}\beta_2 x; q)_{\infty}}{}_3\phi_2 \left( \begin{array}{c} q\beta_2/\alpha_2, q^{\lambda}, \alpha_1 x \\ q^{\lambda+1}\alpha_1/\beta_1, \beta_2 x \end{array} ;q, \alpha_2/\beta_1 \right) \label{eq:fn433-1},\\
&f_{34}(x) = x^{\lambda}\frac{(\beta_2 x;q)_{\infty}}{(q^{-\lambda}\beta_2 x; q)_{\infty}}{}_3\phi_2 \left( \begin{array}{c} \beta_2/\alpha_1, q\beta_2/\alpha_2, q^{\lambda} \\ q\beta_1\beta_2/(\alpha_1\alpha_2), \beta_2 x \end{array} ;q, q \right) \label{eq:fn434-1},\\
&f_{35}(x) = x^{\lambda}\frac{(\alpha_2 x/q;q)_{\infty}}{(q^{-\lambda-1}\beta_1\beta_2 x/\alpha_1; q)_{\infty}}{}_3\phi_2 \left( \begin{array}{c} q^{\lambda}\alpha_1\alpha_2/(\beta_1\beta_2), q\alpha_1/\beta_1, q\alpha_1/\beta_2 \\ q^{\lambda+2}\alpha_1/(\beta_1\beta_2 x), q\alpha_1\alpha_2/(\beta_1\beta_2) \end{array} ;q, q/(\alpha_1 x) \right) \label{eq:fn435-1},\\
&f_{36}(x) = x^{\lambda}\frac{(\alpha_1 x, \alpha_2 x/q;q)_{\infty}}{(\beta_2 x/q, q^{-\lambda}\beta_1 x; q)_{\infty}}{}_3\phi_2 \left( \begin{array}{c} q^{-\lambda+1}, \alpha_2/\beta_2,  q/(\alpha_1 x)\\ q^{-\lambda+1}\beta_1/\alpha_1, q^2/(\beta_2 x) \end{array} ;q, q \right) \label{eq:fn436-1},\\
&f_{37}(x) = x^{\lambda}{}_3\phi_2 \left( \begin{array}{c} q^{\lambda}\alpha_1\alpha_2/(\beta_1\beta_2), q^{\lambda}, \alpha_1 x \\ q^{\lambda+1}\alpha_1/\beta_1, q^{\lambda+1}\alpha_1/\beta_2 \end{array} ;q, q \right) \label{eq:fn437-1}.
\end{align}
By applying formulas for $q$-hypergeometric series to functions (\ref{eq:fn432-1})-(\ref{eq:fn437-1}), we can obtain the relationships with functions (\ref{eq:yb1m0})-(\ref{eq:yxm0}).
We describe them in the following propositions.
\begin{rmk}
Judging from the form of the $q$-hypergeometric series, we cannot rule out the possibility that the functions $f_{34}(x)$, $f_{37}(x)$, and $f_{36}(x)$ are solutions to the non-homogeneous variant of $q$-hypergeometric equation of degree 2.
However, the difference between any two such solutions is a solution to the homogeneous variant of $q$-hypergeometric equation of degree 2 given by (\ref{eq:g1eqdeg2k}).
\end{rmk}
Note that equation (\ref{eq:e2}) is symmetric about $a_1$ and $a_2$, and about $b_1$ and $b_2$.
So the functions replaced $\beta_1$ and $\beta_2$, and $q\alpha_1$ and $\alpha_2$ in above functions are also the solutions.
In this subsection, $\check{f}(x)$ is assumed to represent the function replaced $\beta_1$ and $\beta_2$ in $f(x)$.

We can apply the formula \cite[\rm(I\hspace{-,15em}I\hspace{-,15em}I-10)]{GR}
\begin{align}
{}_3\phi_2 \left( \begin{array}{c} a, b, c\\ d, e \end{array} ;q, \frac{de}{abc} \right) 
= \frac{(b, de/(ab), de/(bc) ;q)_{\infty}}{(d, e, de/(abc) ;q)_{\infty}}{}_3\phi_2 \left( \begin{array}{c} d/b, e/b, de/(abc)\\ de/(ab), de/(bc) \end{array} ;q, b \right) \label{GRIII-10}
\end{align}
to functions (\ref{eq:fn432-1}) and (\ref{eq:fn435-1}), and obtain the following proposition.
\begin{prop} \label{prop:f32}
The functions $f_{32}(x)$ and $\check{f}_{32}(x)$ are constant multiples of $y_{\beta_2}(x)$ and
$y_{\beta_1}(x)$, respectively.
The function $f_{35}(x)$ is a pseudo-constant multiple of $y_x(x)$.
\begin{proof}
Set $a=q^{\lambda}$, $b=q^{\lambda}\alpha_1\alpha_2/(\beta_1\beta_2)$, $c=q^{\lambda+1}/(\beta_1 x)$, $d=q^{\lambda+1}\alpha_1/\beta_1$, $e=q^{\lambda}\alpha_2/\beta_1$ in (\ref{eq:fn432-1}), and apply formula (\ref{GRIII-10}).
Then we have
\begin{equation}
\begin{split}
&f_{32}(x) = S_1 y_{\beta_2}(x),\\
&\qquad S_1 = \frac{\beta_2^{-\lambda+\mu'}}{1-q} \frac{(q^{\lambda}\alpha_1\alpha_2/(\beta_1\beta_2), \beta_2/\alpha_1, q\beta_2/\alpha_2 ;q)_{\infty}}{(q^{\lambda+1}\alpha_1/\beta_1, q^{\lambda}\alpha_2/\beta_1, q ;q)_{\infty}}. \label{eq:f32yb2}
\end{split}
\end{equation}
Replacing $\beta_1$ and $\beta_2$ in (\ref{eq:f32yb2}) yields
\begin{align}
\check{f}_{32}(x) = \check{S}_1 y_{\beta_1}(x).
\label{eq:fc32yb1}
\end{align}
Set $a=q\alpha_1/\beta_1$, $b=q^{\lambda}\alpha_1\alpha_2/(\beta_1\beta_2)$, $c=q\alpha_1/\beta_2$, $d=q^{\lambda+2}\alpha_1/(\beta_1\beta_2 x)$, $e=q\alpha_1\alpha_2/(\beta_1\beta_2)$ in (\ref{eq:fn435-1}).
By formula (\ref{GRIII-10}), we obtain
\begin{equation}
\begin{split}
&f_{35}(x) = S_2(x) y_x(x),  \\
&\qquad S_2(x) = \frac{q^{-\lambda+\mu'}}{1-q} \frac{(q^{\lambda}\alpha_1\alpha_2/(\beta_1\beta_2), q^{-\lambda+1} ;q)_{\infty}\vartheta_q(\alpha_2 x/q)}{(q\alpha_1\alpha_2/(\beta_1\beta_2), q ;q)_{\infty}\vartheta_q(q^{-\lambda-1}\beta_1\beta_2 x/\alpha_1)} x^{\lambda-\mu'}, \label{eq:f35yx} \\
\end{split}
\end{equation}
where $\vartheta_q(t) = (t, q/t, q;q)_{\infty}$.
The function $S_2(x)$ is a pseudo-constant that satisfies $S_2(qx) = S_2(x)$.
\end{proof}
\end{prop}
We can obtain the following propositions by using \cite[formula (3.3.3)]{GR}
\begin{align}
&_3\phi_2 \left( \begin{array}{c} a,b,c \\ d,e \end{array} ; q, \frac{de}{abc} \right) 
=\frac{(e/b,e/c,cq/a,q/d;q)_{\infty}}{(e,cq/d,q/a,e/(bc);q)_{\infty}}
{}_3\phi_2 \left( \begin{array}{c} c,d/a,cq/e \\ cq/a,bcq/e \end{array} ; q, \frac{bq}{d} \right) \nonumber \\
&\qquad -\frac{(q/d,eq/d,b,c,d/a,de/(bcq),bcq^2/(de);q)_{\infty}}{(d/q,e,bq/d,cq/d,q/a,e/(bc),bcq/e;q)_{\infty}}
{}_3\phi_2 \left( \begin{array}{c} aq/d,bq/d,cq/d \\ q^2/d,eq/d \end{array} ; q, \frac{de}{abc} \right), \label{GR333}
\end{align}
where $|bq/d|<1$ and $|de/(abc)|<1$,
and \cite[formula (3.3.1)]{GR}
\begin{align}
&_3\phi_2 \left( \begin{array}{c} a,b,c \\ d,e \end{array} ;q,\frac{de}{abc} \right)
=\frac{(e/b,e/c;q)_{\infty}}{(e,e/(bc);q)_{\infty}} {}_3\phi_2 \left( \begin{array}{c} d/a,b,c \\ d,bcq/e \end{array} ;q,q \right) \notag \\
&\hspace{50pt} +\frac{(d/a,b,c,de/(bc);q)_{\infty}}{(d,e,bc/e,de/(abc);q)_{\infty}} {}_3\phi_2 \left( \begin{array}{c} e/b,e/c,de/(abc) \\ de/(bc),eq/(bc) \end{array} ;q,q \right), \label{GR331}
\end{align}
where $|de/(abc)|<1$.
\begin{prop} \label{prop:f33}
The following relation holds:
\begin{equation}
\begin{split}
f_{33}(x) &= \frac{\check{S}_1 + A_{1,2}K_{6,3}}{A_{1,1}} y_{\beta_1}(x) + \frac{A_{1,2}}{A_{1,1}} y_{\beta_2}(x), \\
&A_{1,1} = \frac{(\beta_1/\alpha_1, \alpha_2/\beta_2, q^{-\lambda}\beta_2/\alpha_1; q)_{\infty}}{(q^{\lambda}\alpha_2/\beta_2, \beta_2/\alpha_1, q^{-\lambda}\beta_1/\alpha_1; q)_{\infty}}, \\
&A_{1,2} = -\frac{\beta_2^{-\lambda+\mu'}}{1-q} \frac{(q^{-\lambda}\beta_2/\alpha_1, q^{\lambda}\alpha_1\alpha_2/(\beta_1\beta_2), q\beta_2/\alpha_2, \beta_1/\alpha_1, \alpha_2/\beta_2; q)_{\infty}}{(q^{\lambda}\alpha_2/\beta_2, \alpha_2/\beta_1, q^{-\lambda}\beta_1/\alpha_1, q^{\lambda+1}\alpha_1/\beta_1, q; q)_{\infty}}, \\
&K_{6,3}= -\left( \frac{\beta_2}{\beta_1} \right)^{\lambda-\mu'}
\frac{\vartheta_q(\beta_1/\beta_2)\vartheta_q(\alpha_2/\beta_1)\vartheta_q(q^{\lambda}\alpha_1/\beta_1)}{\vartheta_q(\beta_2/\beta_1)\vartheta_q(\alpha_2/\beta_2)\vartheta_q(q^{\lambda}\alpha_1/\beta_2)}. 
\label{eq:f33yb1yb2}
\end{split}
\end{equation}
\begin{proof}
By applying formula (\ref{GR333}) to $\check{f}_{32}(x)$, obtained by replacing $\beta_1$ and $\beta_2$ in (\ref{eq:fn432-1}), and to (\ref{eq:yb2m0}), we have
\begin{equation}
\check{f}_{32}(x) = A_{1,1}f_{33}(x) + A_{1,2} y_2^{(10)}(x),
\label{eq:fc32-f33} 
\end{equation}  
where $a = q^{\lambda+1}/(\beta_2 x), b = q^{\lambda}\alpha_1\alpha_2/(\beta_1\beta_2), c = q^{\lambda}, d = q^{\lambda+1}\alpha_1/\beta_1, e = q^{\lambda}\alpha_2/\beta_1$, and
\begin{equation}
y_{\beta_2}(x) = y_2^{(10)}(x) + K_{6,3}y_{\beta_1}(x) \label{eq:yb2-y2(10)},
\end{equation}
where $a = q\beta_2/\alpha_2, b = q^{-\lambda}\beta_2 x, c = \beta_2/\alpha_1, d = q\beta_2/\beta_1, e = \beta_2 x$. 
For an explicit formula of $y_2^{(10)}(x)$, see Appendix \ref{app:y1y2m=0}.
Equations (\ref{eq:fc32yb1}), (\ref{eq:fc32-f33}), and (\ref{eq:yb2-y2(10)}) lead to equation (\ref{eq:f33yb1yb2}).
\end{proof}
\end{prop}
\begin{prop} \label{prop:f343736}
The following relations hold:
\begin{equation}
\begin{split}
&f_{34}(x) - \frac{C_1(x)}{S_3} = \frac{\check{S}_1 + A_{1,2}K_{6,3}}{A_{1,1}S_3} y_{\beta_1}(x) - \frac{A_{1,2}}{A_{1,1}S_3} y_{\beta_2}(x) , \\
&\qquad S_3 = -\frac{(q^{\lambda}\alpha_1\alpha_2/(\beta_1\beta_2), q\alpha_1/\beta_1; q)_{\infty}}{(q^{\lambda+1}\alpha_1/\beta_1, \alpha_1\alpha_2/(\beta_1\beta_2); q)_{\infty}}, \\
&\qquad C_1(x) = x^{\lambda}\frac{(\beta_2/\alpha_1, q\beta_2/\alpha_2, q^{\lambda}, \alpha_1\alpha_2 x/\beta_1 ;q)_{\infty}}{(q^{\lambda+1}\alpha_1/\beta_1, \beta_1\beta_2/(\alpha_1\alpha_2), \alpha_2/\beta_1, q^{-\lambda}\beta_2 x  ;q)_{\infty}}  \\
&\hspace{120pt} \cdot {}_3\phi_2 \left( \begin{array}{c} q^{\lambda}\alpha_1\alpha_2/(\beta_1\beta_2), q\alpha_1/\beta_1, \alpha_2/\beta_1 \\ \alpha_1\alpha_2 x/\beta_1, q\alpha_1\alpha_2/(\beta_1\beta_2) \end{array} ;q, q \right),
\label{rl:f34}
\end{split}
\end{equation}
\begin{equation}
\begin{split}
&f_{37}(x) - \frac{C_2(x)}{S_4} = \frac{S_1}{S_4} y_{\beta_2}(x) , \\
&\qquad S_4 = -\frac{(\beta_2/\alpha_1, \alpha_2/\beta_1; q)_{\infty}}{(q^{\lambda}\alpha_2/\beta_1, q^{-\lambda}\beta_2/\alpha_1; q)_{\infty}}, \\
&\qquad C_2(x) = x^{\lambda}\frac{(\alpha_1 x, q^{\lambda}\alpha_1\alpha_2/(\beta_1\beta_2), q^{\lambda}, q\beta_2/\beta_1; q)_{\infty}}{(q^{\lambda+1}\alpha_1/\beta_1, q^{\lambda}\alpha_2/\beta_1, q^{\lambda}\alpha_1/\beta_2, q^{-\lambda}\beta_2 x, \beta_2/\alpha_1, \alpha_2/\beta_1; q)_{\infty}} \\
&\hspace{120pt}\cdot{}_3\phi_2 \left( \begin{array}{c} \beta_2/\alpha_1, \alpha_2/\beta_1, q^{-\lambda}\beta_2 x \\ q\beta_2/\beta_1, q^{-\lambda+1}\beta_2/\alpha_1 \end{array} ;q, q \right),\label{rl:f37}
\end{split}
\end{equation}
\begin{equation}
\begin{split}  
&f_{36}(x) - \frac{C_3(x)}{S_5(x)}= \frac{S_2(x)}{S_5(x)}y_x(x) , \\
&\qquad S_5(x) = -\frac{(q^{\lambda}\alpha_1\alpha_2/(\beta_1\beta_2), q\alpha_1/\beta_1 ;q)_{\infty} \vartheta_q(q^{-\lambda}\beta_1 x)\vartheta_q(\beta_2 x/q)}{(q\alpha_1\alpha_2/(\beta_1\beta_2), q^{\lambda}\alpha_1/\beta_1 ;q)_{\infty}\vartheta_q(q^{-\lambda-1}\beta_1\beta_2 x/\alpha_1)\vartheta_q(\alpha_1 x)}, \\
&\qquad C_3(x) = x^{\lambda}\frac{(\alpha_2 x/q, q^{-\lambda+1}, \alpha_2/\beta_2 ;q)_{\infty}}{(q^{-\lambda-1}\beta_1\beta_2 x/\alpha_1, q\alpha_1\alpha_2/(\beta_1\beta_2), q^{-\lambda}\beta_1/\alpha_1 ;q)_{\infty}} \\
&\hspace{120pt} \cdot {}_3\phi_2 \left( \begin{array}{c} q^{\lambda+1}/(\beta_1 x), q^{\lambda}\alpha_1\alpha_2/(\beta_1\beta_2), q\alpha_1/\beta_1 \\ q^{\lambda+2}\alpha_1/(\beta_1\beta_2 x), q^{\lambda+1}\alpha_1/\beta_1 \end{array} ;q, q \right).
\label{rl:f36}
\end{split}
\end{equation}
\begin{proof}
By applying formula (\ref{GR331}) to finctions (\ref{eq:fn433-1}), (\ref{eq:fn432-1}), and (\ref{eq:fn435-1}), we have
\begin{equation}
f_{33}(x) = S_3 f_{34}(x) + C_1(x) , \label{eq:f33-331} 
\end{equation}
where $a = \alpha_1 x, b = q\beta_2/\alpha_2, c = q^{\lambda}, d = \beta_2 x, e = q^{\lambda+1}\alpha_1/\beta_1$ ,
\begin{equation}
f_{32}(x) = S_4 f_{37}(x) + C_2(x) , \label{eq:f32-331} 
\end{equation}
where $a = q^{\lambda+1}/(\beta_1 x), b = q^{\lambda}\alpha_1\alpha_2/(\beta_1\beta_2), c = q^{\lambda}, d = q^{\lambda+1}\alpha_1/\beta_1, e = q^{\lambda}\alpha_2/\beta_1$ , and 
\begin{equation}
f_{35}(x) = C_3(x) + S_5(x)f_{36}(x) , \label{eq:f35-331}
\end{equation}
where $a = q\alpha_1/\beta_2, b = q^{\lambda}\alpha_1\alpha_2/(\beta_1\beta_2), c = q\alpha_1/\beta_1, d = q^{\lambda+2}\alpha_1/(\beta_1\beta_2 x), e = q\alpha_1\alpha_2/(\beta_1\beta_2)$.
We obtain (\ref{rl:f34}) from (\ref{eq:f33yb1yb2}) and (\ref{eq:f33-331}), (\ref{rl:f37}) from (\ref{eq:f32yb2}) and (\ref{eq:f32-331}), and (\ref{rl:f36}) from (\ref{eq:f35yx}) and (\ref{eq:f35-331}).
\end{proof}
\end{prop}
If $f_{34}(x)$, $f_{37}(x)$, and $f_{36}(x)$ satisfy the non-homogeneous variant of $q$-hypergeometric equation of degree 2, then $C_1(x)/S_3$, $C_2(x)/S_4$, and $C_3(x)/S_5(x)$ also satisfy the same equation.

\section{Linear relationships among integral solutions associated with $q$-middle convolution}$ $ \label{sec:LinRel} 

To the $_3\phi_2$ series in functions (\ref{eq:yb1m0}), (\ref{eq:yb2m0}), and (\ref{eq:yxm0}), we can apply formula (\ref{GR333}).
We transform $_3\phi_2(q^{-\lambda}\beta_1 x, \beta_1/\alpha_1,q\beta_1/\alpha_2 ; \beta_1 x, q\beta_1/\beta_2 ; q,  
q^{\lambda}\alpha_1\alpha_2/(\beta_1\beta_2))$ in function (\ref{eq:yb1m0}).
Set $a=q^{-\lambda}\beta_1 x$, $b=\beta_1/\alpha_1$, $c=q\beta_1/\alpha_2$, $d=\beta_1 x$, and $e=q\beta_1/\beta_2$.
Then we have
\begin{align}
&y_{\beta_1}(x) \notag \\
&\quad = (1-q)\beta_1^{\lambda-\mu'}x^{\lambda}\frac{(\beta_1 x, q, q\alpha_1/\beta_2, \alpha_2/\beta_2, q^{\lambda+2}/(\alpha_2 x), q/(\beta_1 x) ;q)_{\infty}}{(q^{-\lambda}\beta_1 x, \beta_1/\alpha
_1, q\beta_1/\alpha_2, q^2/(\alpha_2 x), q^{\lambda+1}/(\beta_1 x), \alpha_1\alpha_2/(\beta_1\beta_2)  ;q)_{\infty}} \nonumber \\ 
&\hspace{150pt}\cdot\!_3\phi_2 \left( \begin{array}{c} q\beta_1/\alpha_2,q^{\lambda},q\beta_2/\alpha_2 \\ q^{\lambda+2}/(\alpha_2 x), q\beta_1\beta_2/(\alpha_1\alpha_2)  \end{array} ;q, q/(\alpha_1 x) \right) \nonumber \\
&\quad -(1-q)\beta_1^{\lambda-\mu'}x^{\lambda} \nonumber \\
&\qquad\cdot \frac{(\beta_1 x, q, q/(\beta_1 x), q^2/(\beta_2 x), q^{\lambda}, \alpha_1\alpha_2 x/(q\beta_2), q^2\beta_2/(\alpha_1\alpha_2 x) ;q)_{\infty}}{(q^{-\lambda}\beta_1 x, \beta_1 x/q, q/(\alpha_1 x), q^2/(\alpha_2 x), q^{\lambda+1}/(\beta_1 x), \alpha_1\alpha_2/(\beta_1\beta_2), q\beta_1\beta_2/(\alpha_1\alpha_2)   ;q)_{\infty}} \nonumber \\
&\hspace{150pt}\cdot\!_3\phi_2 \left( \begin{array}{c} q^{-\lambda+1},q/(\alpha_1 x),q^2/(\alpha_2 x) \\ q^2/(\beta_1 x), q^2/(\beta_2 x)  \end{array} ;q, q^{\lambda}\frac{\alpha_1\alpha_2}{\beta_1\beta_2} \right) \notag \\ 
&\quad = (1-q)\beta_1^{\lambda-\mu'}x^{\lambda}\frac{(\beta_1 x, q, q\alpha_1/\beta_2, \alpha_2/\beta_2, q^{\lambda+2}/(\alpha_2 x), q/(\beta_1 x) ;q)_{\infty}}{(q^{-\lambda}\beta_1 x, \beta_1/\alpha
_1, q\beta_1/\alpha_2, q^2/(\alpha_2 x), q^{\lambda+1}/(\beta_1 x), \alpha_1\alpha_2/(\beta_1\beta_2)  ;q)_{\infty}} \nonumber \\ 
&\hspace{150pt}\cdot\!_3\phi_2 \left( \begin{array}{c} q\beta_1/\alpha_2,q^{\lambda},q\beta_2/\alpha_2 \\ q^{\lambda+2}/(\alpha_2 x), q\beta_1\beta_2/(\alpha_1\alpha_2)  \end{array} ;q, q/(\alpha_1 x) \right) \nonumber \\
&\quad -q^{-\lambda+\mu'}\beta_1^{\lambda-\mu'}x^{\lambda-\mu'}
\frac{\vartheta_q(\beta_1 x)\vartheta_q(q^{\lambda})\vartheta_q(\alpha_1\alpha_2 x/(q\beta_2))}{\vartheta_q(q^{-\lambda}\beta_1 x)\vartheta_q(\beta_1 x/q)\vartheta_q(\alpha_1\alpha_2/(\beta_1\beta_2))} \; y_x(x).
\label{eq:b1333} 
\end{align}

\begin{rmk} \label{rem:b1b2xm=0}
Let $J_{1,3}(x)$ be the coefficient of $y_x(x)$ in the second term on the right-hand side of (\ref{eq:b1333}).
Since the coefficient $J_{1,3}(x)$ is a pseudo-constant which satisfies $J_{1,3}(qx) = J_{1,3}(x)$, $J_{1,3}(x)y_x(x)$ is also a solution to equation (\ref{eq:g1eqdeg2k}).
This implies that the first term on the right-hand side is also a solution to equation (\ref{eq:g1eqdeg2k}).
By varying the parameter settings, we can derive 12 solutions that have different expressions from (\ref{eq:yb1m0})-(\ref{eq:yxm0}).
Considering transformations for $y_{\beta_2}(x)$ and $y_x(x)$, we obtain 36 solutions in total with novel expressions.
Excluding duplicates, mere pseudo-constant multiples, and those transferred to other solutions by transformation, 6 new solutions remain (see Appendix \ref{app:y1y2m=0}).
\end{rmk}

To the series $_3\phi_2( q\beta_1/\alpha_2,q^{\lambda},q\beta_2/\alpha_2 ; q^{\lambda+2}/(\alpha_2 x), q\beta_1\beta_2/(\alpha_1\alpha_2) ;q, q/(\alpha_1 x))$ in the first term on the right-hand side of equation (\ref{eq:b1333}), we can apply formula (\ref{GR331}).
Set $a=q^{\lambda}$, $b=q\beta_1/\alpha_2$, $c=q\beta_2/\alpha_2$, $d=q^{\lambda+2}/(\alpha_2 x)$, and $e=q\beta_1\beta_2/(\alpha_1\alpha_2)$.
Then the first term on the right-hand side of equation (\ref{eq:b1333}) is converted as follows:
\begin{align}
&\quad (1-q)\beta_1^{\lambda-\mu'}x^{\lambda} \notag \\
&\qquad\cdot \frac{(\beta_1 x, q, q\alpha_1/\beta_2, \alpha_2/\beta_2, q^{\lambda+2}/(\alpha_2 x), q/(\beta_1 x), \beta_2/\alpha_1 ;q)_{\infty}}{(q^{-\lambda}\beta_1 x, q\beta_1/\alpha_2, q^2/(\alpha_2 x), q^{\lambda+1}/(\beta_1 x), \alpha_1\alpha_2/(\beta_1\beta_2), q\beta_1\beta_2/(\alpha_1\alpha_2), \alpha_2/(q\alpha_1) ;q)_{\infty}} \nonumber \\
&\hspace{150pt}\cdot\!_3\phi_2 \left( \begin{array}{c} q^2/(\alpha_2 x),q\beta_1/\alpha_2,q\beta_2/\alpha_2 \\ q^{\lambda+2}/(\alpha_2 x), q^2\alpha_1/\alpha_2  \end{array} ;q, q \right) \nonumber \\
&\quad +(1-q)\beta_1^{\lambda-\mu'}x^{\lambda}  \nonumber \\
&\qquad\cdot \frac{(\beta_1 x, q, q\alpha_1/\beta_2, \alpha_2/\beta_2, q/(\beta_1 x), q\beta_2/\alpha_2, q^{\lambda+1}/(\alpha_1 x) ;q)_{\infty}}{(q^{-\lambda}\beta_1 x, \beta_1/\alpha
_1, q^{\lambda+1}/(\beta_1 x), \alpha_1\alpha_2/(\beta_1\beta_2), q\beta_1\beta_2/(\alpha_1\alpha_2), q\alpha_1/\alpha_2, q/(\alpha_1 x) ;q)_{\infty}} \nonumber  \\
&\hspace{150pt}\cdot\!_3\phi_2 \left( \begin{array}{c} \beta_2/\alpha_1,\beta_1/\alpha_1,q/(\alpha_1 x) \\ q^{\lambda+1}/(\alpha_1 x), \alpha_2/\alpha_1  \end{array} ;q, q \right) \nonumber \\
&= -\frac{\alpha_2}{q\alpha_1}\beta_1^{\lambda-\mu'}x^{\lambda}
\frac{\vartheta_q(\beta_1 x)\vartheta_q(\alpha_2/\beta_2)\vartheta_q(\beta_2/\alpha_1)}{\vartheta_q(q^{-\lambda}\beta_1 x)\vartheta_q(\alpha_1\alpha_2/(\beta_1\beta_2))\vartheta_q(\alpha_2/(q\alpha_1))} \; y_{\alpha_2}(x) \notag  \\
&\qquad\qquad -\beta_1^{\lambda-\mu'}x^{\lambda}
\frac{\vartheta_q(\beta_1 x)\vartheta_q(\alpha_2/\beta_2)\vartheta_q(\beta_2/\alpha_1)}{\vartheta_q(q^{-\lambda}\beta_1 x)\vartheta_q(\alpha_1\alpha_2/(\beta_1\beta_2))\vartheta_q(\alpha_2/\alpha_1)}\; y_{\alpha_1}(x).
\label{eq:ya2+ya1}
\end{align}
Let $J_{1,1}(x)$ and $J_{1,2}(x)$ be the coefficients of $y_{\alpha_1}(x)$ and $y_{\alpha_2}(x)$ in equation (\ref{eq:ya2+ya1}), respectively.
These coefficients are pseudo-constants.
Summing up the results of calculations so far, we can write
\begin{equation*}
y_{\beta_1}(x) = J_{1,1}(x)y_{\alpha_1}(x) + J_{1,2}(x)y_{\alpha_2}(x) + J_{1,3}(x)y_x(x),
\end{equation*}
where $J_{1,1}(x)$, $J_{1,2}(x)$, and $J_{1,3}(x)$ are pseudo-constants. 
Since $y_{\beta_1}(x)$ and $J_{1,3}(x)y_x(x)$ satisfy the homogeneous equation given by (\ref{eq:g1eqdeg2k}) , and $J_{1,1}(x)y_{\alpha_1}(x)$ and $J_{1,2}(x)y_{\alpha_2}(x)$ satisfy the non-homogeneous equation given by (\ref{eq:y1nonhom}), the coefficients $J_{1,1}(x)$ and $J_{1,2}(x)$ should satisfy the equation
\begin{align*}
J_{1,1}(x)+J_{1,2}(x)=0. 
\end{align*}
It can be verified by a simple calculation.
In conclusion we have
\begin{equation*}
y_{\beta_1}(x) = J_{1,1}(x) (y_{\alpha_1}(x) - y_{\alpha_2}(x)) + J_{1,3}(x)y_x(x).
\end{equation*}

By executing the identical procedure of calculations on $y_{\beta_1}(x)$, $y_{\beta_2}(x)$, and $y_x(x)$ with variations in parameterization, 9 relational expressions are derived.
\begin{thm} \label{thm:rlm=0}
For $y_{\alpha_1}(x)$, $y_{\alpha_2}(x)$, $y_{\lambda}(x)$, $y_{\beta_1}(x)$, $y_{\beta_2}(x)$, and $y_x(x)$, the linear relationships
\begin{align*}
y_{\beta_1}(x) &= J_{1,1}(x) (y_{\alpha_1}(x) - y_{\alpha_2}(x)) + J_{1,3}(x)y_x(x) \\ 
               &= J_{2,1}(x) (y_{\alpha_1}(x) - y_{\alpha_2}(x)) + J_{2,3}(x)y_{\beta_2}(x)\\
               &= J_{3,1}(x) (y_{\alpha_2}(x) - y_{\lambda}(x))+ J_{3,3}(x)y_x(x)\\
               &= J_{4,1}(x) (y_{\alpha_2}(x) - y_{\lambda}(x)) + J_{4,3}(x)y_{\beta_2}(x)\\
               &= J_{5,1}(x) (y_{\alpha_1}(x) - y_{\lambda}(x)) + J_{5,3}(x)y_x(x)\\
               &= J_{6,1}(x) (y_{\alpha_1}(x) - y_{\lambda}(x)) + J_{6,3}(x)y_{\beta_2}(x),
\end{align*}
\begin{align*}                            
y_{\beta_2}(x) &= K_{1,1}(x) (y_{\alpha_1}(x) - y_{\alpha_2}(x)) + K_{1,3}(x)y_x(x)\\
               &= K_{3,1}(x) (y_{\alpha_2}(x)  - y_{\lambda}(x)) + K_{3,3}(x)y_x(x)\\
               &= K_{5,1}(x) ( y_{\alpha_1}(x) - y_{\lambda}(x) ) + K_{5,3}(x)y_x(x)
\end{align*}
hold, where $J_{m,n}(x)$\ $(m=1,2,3,4,5,6,\ n=1,3)$ and $K_{m',n'}(x)$\ $(m'=1,3,5,\ n'=1,3)$ are pseudo-constants (see Appendix \ref{coef:JK}).
\end{thm}
The relation among $y_{\beta_1}(x)$, $y_{\beta_2}(x)$, and $y_x(x)$ follows directly from the first two equations of Theorem \ref{thm:rlm=0}.
\begin{cor}
The relation
\begin{align}
y_x(x) =\frac{J_{2,1}(x) - J_{1,1}(x)}{J_{1,3}(x)J_{2,1}(x)} y_{\beta_1}(x) + \frac{J_{1,1}(x)J_{2,3}(x)}{J_{1,3}(x)J_{2,1}(x)} y_{\beta_2}(x) \label{eq:yx=yb1+yb2}
\end{align}
holds for $y_{\beta_1}(x)$, $y_{\beta_2}(x)$, and $y_x(x)$ which satisfy equation (\ref{eq:g1eqdeg2k}).
\end{cor}

\section{Concluding remarks}
In this paper, we have shown correspondences between integral solutions obtained by applying $q$-middle convolution and solutions obtained by alternative methods.
The presence of correspondence between solutions obtained through different methods implies the validity of each respective method.
And from current results, we can conclude that $q$-middle convolution can be used as a tool to analyze $q$-difference equations and their solutions.

In \cite{AT}, for the first-order $q$-difference equation satisfied by the function $y(x) = x^{\mu}(\alpha_1 x, \\
\alpha_2 x; q)_{\infty}/(\beta_1 x, \beta_2 x;q)_{\infty}$, we applied the $q$-middle convolution under the condition $q^{\lambda}=q^{\mu}\beta_1\beta_2/(\alpha_1\alpha_2)$.
The resulting single second-order $q$-difference equations correspond to the variant of $q$-hypergeometric equation of degree 2.
As in Section \ref{sec:is-vqhg2}, their integral solutions can be expressed as $q$-hypergeomtric series ${}_3\phi_2$ by specializing the value of the parameter $\xi$ (see \cite[Section 4.1.2]{AT}).
We conducted the same investigations as in Section \ref{sec:correso} and Section \ref{sec:LinRel} for these integral solutions and obtained similar results.
Here, we omit these results.
%

In \cite{AT}, integral solutions to the variant of $q$-hypergeometric equation of degree 3 were also obtained by the $q$-middle convolution.
These solutions are expressed by using the $q$-hypergeometric series ${}_4\phi_3$.
They should also be investigated, but since no suitable formula for ${}_4\phi_3$ can be found, a different approach is needed.

There is still much more to explore about solutions to the variants of $q$-hypergeometic equation.

\section*{Acknowledgements}
The author would like to express sincere gratitude to Professor Kouichi Takemura for invaluable advices and continuous support.
The author also wishes to thank the reviewers for their constructive comments and helpful suggestions, which significantly contributed to improving the quality of the manuscript.

\appendix

\section{Supplement to Remark \ref{rem:b1b2xm=0}} \label{app:y1y2m=0}
The following functions are the solutions to equation $(\ref{eq:g1eqdeg2k})$ which have different expressions from $y_{\beta_1}$, $y_{\beta_2}$, and $y_x(x)$.
\begin{align*}
&y_1^{(1)}(x) = (1-q)\beta_1^{\lambda-\mu'}x^{\lambda}\frac{(\beta_1 x, q, q\alpha_1/\beta_2, \alpha_2/\beta_2, q^{\lambda+2}/(\alpha_2 x), q/(\beta_1 x) ;q)_{\infty}}{(q^{-\lambda}\beta_1 x, \beta_1/\alpha_1, q\beta_1/\alpha_2, q^2/(\alpha_2 x), q^{\lambda+1}/(\beta_1 x), \alpha_1\alpha_2/(\beta_1\beta_2)  ;q)_{\infty}} \notag \\ 
&\hspace{100pt}\cdot\!_3\phi_2 \left( \begin{array}{c} q\beta_1/\alpha_2,q^{\lambda},q\beta_2/\alpha_2 \\ q^{\lambda+2}/(\alpha_2 x), q\beta_1\beta_2/(\alpha_1\alpha_2)  \end{array} ;q, q/(\alpha_1 x) \right),\quad |q/(\alpha_1 x)|<1, 
\\
&y_1^{(5)}(x) = (1-q)\beta_1^{\lambda-\mu'}x^{\lambda}\frac{(\beta_1 x, q, q^{\lambda+1}/(\beta_2 x), \alpha_2/\beta_2, q^2\alpha_1/\alpha_2, q/(\beta_1 x) ;q)_{\infty}}{(q^{-\lambda}\beta_1 x, \beta_1/\alpha_1, q\beta_1/\alpha_2, q^2/(\alpha_2 x), q\alpha_1/\beta_1, q^{\lambda}\alpha_2/(\beta_1\beta_2 x) ;q)_{\infty}} \notag \\ 
&\hspace{100pt}\cdot\!_3\phi_2 \left( \begin{array}{c} q\beta_1/\alpha_2, \alpha_1 x, q\beta_2/\alpha_2 \\ q^2\alpha_1/\alpha_2, q^{-\lambda+1}\beta_1\beta_2 x/\alpha_2  \end{array} ;q, q^{-\lambda+1} \right),\quad \lambda<1, \\
&y_1^{(9)}(x) = (1-q)\beta_1^{\lambda-\mu'}x^{\lambda}\frac{(\beta_1 x, q, q^{\lambda+1}/(\beta_2 x), q\alpha_1/\beta_2, \alpha_2/\alpha_1, q/(\beta_1 x) ;q)_{\infty}}{(q^{-\lambda}\beta_1 x, \beta_1/\alpha_1, q\beta_1/\alpha_2, q/(\alpha_1 x), \alpha_2/\beta_1, q^{\lambda+1}\alpha_1/(\beta_1\beta_2 x) ;q)_{\infty}} \notag \\ 
&\hspace{100pt}\cdot\!_3\phi_2 \left( \begin{array}{c} \beta_1/\alpha_1, \alpha_2 x/q, \beta_2/\alpha_1 \\ \alpha_2/\alpha_1, q^{-\lambda}\beta_1\beta_2 x/\alpha_1  \end{array} ;q, q^{-\lambda+1} \right),\quad \lambda<1, \\
&y_2^{(2)}(x) = (1-q)\beta_2^{\lambda-\mu'}x^{\lambda}\frac{(q, q\beta_2/\beta_1, \alpha_1 x, \alpha_2 x/q, q^{\lambda+2}/(\alpha_2 x), \beta_1/\beta_2 ;q)_{\infty}}{(q^{-\lambda}\beta_2 x, \beta_2/\alpha_1, q\beta_2/\alpha_2, q\beta_1/\alpha_2, q^{\lambda+1}/(\beta_2 x), \alpha_1\alpha_2 x/(q\beta_2)  ;q)_{\infty}} \notag \\ 
&\hspace{100pt}\cdot\!_3\phi_2 \left( \begin{array}{c} q\beta_2/\alpha_2,q^{\lambda+1}/(\beta_1 x), q^2/(\alpha_2 x) \\ q^{\lambda+2}/(\alpha_2 x), q^2\beta_2/(\alpha_1\alpha_2 x)  \end{array} ;q, \beta_1/\alpha_1  \right),\quad |\beta_1/\alpha_1|<1, \\
&y_2^{(6)}(x) = (1-q)\beta_2^{\lambda-\mu'}x^{\lambda}\frac{(q, q\beta_2/\beta_1, q^{\lambda}, \alpha_2 x/q, q^2\alpha_1/\alpha_2, \beta_1/\beta_2 ;q)_{\infty}}{(q^{-\lambda}\beta_2 x, \beta_2/\alpha_1, q\beta_2/\alpha_2, q\beta_1/\alpha_2, q\alpha_1/\beta_2, q^{\lambda-1}\alpha_2/\beta_2 ;q)_{\infty}} \notag \\ 
&\hspace{100pt}\cdot\!_3\phi_2 \left( \begin{array}{c} q\beta_2/\alpha_2, q\alpha_1/\beta_1, q^2/(\alpha_2 x) \\ q^2\alpha_1/\alpha_2, q^{-\lambda+2}\beta_2/\alpha_2  \end{array} ;q, q^{-\lambda}\beta_1 x \right),\quad |q^{-\lambda}\beta_1 x|<1, \\
&y_2^{(10)}(x) = (1-q)\beta_2^{\lambda-\mu'}x^{\lambda}\frac{(q, q\beta_2/\beta_1, q^{\lambda}, \alpha_1 x, \alpha_2/\alpha_1, \beta_1/\beta_2 ;q)_{\infty}}{(q^{-\lambda}\beta_2 x, \beta_2/\alpha_1, q\beta_2/\alpha_2, \beta_1/\alpha_1, \alpha_2/\beta_2, q^{\lambda}\alpha_1/\beta_2 ;q)_{\infty}} \notag \\ 
&\hspace{100pt}\cdot\!_3\phi_2 \left( \begin{array}{c} \beta_2/\alpha_1, \alpha_2/\beta_1, q/(\alpha_1 x) \\ \alpha_2/\alpha_1, q^{-\lambda+1}\beta_2/\alpha_1  \end{array} ;q, q^{-\lambda}\beta_1 x \right),\quad |q^{-\lambda}\beta_1 x|<1.
\end{align*}

\section{Supplement to Theorem \ref{thm:rlm=0}} \label{coef:JK}
The following functions are the coefficients of the linear relations shown in Theorem \ref{thm:rlm=0}.
\begin{align*}
&J_{1,1}(x) = -\beta_1^{\lambda-\mu'}x^{\lambda}
\frac{\vartheta_q(\beta_1 x)\vartheta_q(\alpha_2/\beta_2)\vartheta_q(\beta_2/\alpha_1)}{\vartheta_q(q^{-\lambda}\beta_1 x)\vartheta_q(\alpha_1\alpha_2/(\beta_1\beta_2))\vartheta_q(\alpha_2/\alpha_1)},\\
&J_{1,3}(x)= -q^{-\lambda+\mu'}\beta_1^{\lambda-\mu'}x^{\lambda-\mu'}
\frac{\vartheta_q(\beta_1 x)\vartheta_q(q^{\lambda})\vartheta_q(\alpha_1\alpha_2 x/(q\beta_2))}{\vartheta_q(q^{-\lambda}\beta_1 x)\vartheta_q(\beta_1 x/q)\vartheta_q(\alpha_1\alpha_2/(\beta_1\beta_2))}, 
\end{align*}
\begin{align*}
&J_{2,1}(x) = -\beta_1^{\lambda-\mu'}x^{\lambda}
\frac{\vartheta_q(\alpha_1 x)\vartheta_q(\alpha_2 x/q)\vartheta_q(\beta_2/\beta_1)}{\vartheta_q(q^{-\lambda}\beta_1 x)\vartheta_q(\alpha_1\alpha_2 x/(q\beta_1))\vartheta_q(\alpha_2/\alpha_1)},\\
&J_{2,3}(x)= -\left( \frac{\beta_1}{\beta_2} \right)^{\lambda-\mu'}
\frac{\vartheta_q(\beta_2/\beta_1)\vartheta_q(\alpha_1\alpha_2 x/(q\beta_2))\vartheta_q(q^{-\lambda}\beta_2 x)}{\vartheta_q(q^{-\lambda}\beta_1 x)\vartheta_q(\beta_1 /\beta_2)\vartheta_q(\alpha_1\alpha_2 x/(q\beta_1))}, \\
&J_{3,1}(x) = -\frac{\alpha_2}{q\alpha_1}\beta_1^{\lambda-\mu'}x^{\lambda}
\frac{\vartheta_q(\beta_1 x)\vartheta_q(\alpha_2/\beta_2)\vartheta_q(q^{-\lambda}\beta_2 x)}{\vartheta_q(\beta_1/\alpha_1)\vartheta_q(q^{-\lambda+1}\beta_1\beta_2 x/\alpha_2)\vartheta_q(q^{-\lambda-1}\alpha_2 x)},\\
&J_{3,3}(x)= - \left( \frac{\beta_1}{q} \right)^{\lambda-\mu'} x^{\lambda-\mu'}
\frac{\vartheta_q(\beta_1 x)\vartheta_q(\alpha_1 x)\vartheta_q(q^{\lambda-1}\alpha_1/\beta_2)}{\vartheta_q(\beta_1/\alpha_1)\vartheta_q(\beta_1 x/q)\vartheta_q(q^{-\lambda+1}\beta_1\beta_2 x/\alpha_2)}, \\
&J_{4,1}(x) = -\frac{\alpha_2}{q\alpha_1}\beta_1^{\lambda-\mu'}x^{\lambda}
\frac{\vartheta_q(q^{\lambda})\vartheta_q(\alpha_2 x/q)\vartheta_q(\beta_2/\beta_1)}{\vartheta_q(\beta_1/\alpha_1)\vartheta_q(q^{\lambda-1}\alpha_2/\beta_1)\vartheta_q(q^{-\lambda-1}\alpha_2 x)},\\
&J_{4,3}(x)= -\left( \frac{\beta_1}{\beta_2} \right)^{\lambda-\mu'}
\frac{\vartheta_q(\beta_2/\beta_1)\vartheta_q(q^{\lambda-1}\alpha_2/\beta_2)\vartheta_q(\beta_2/\alpha_1)}{\vartheta_q(\beta_1/\alpha_1)\vartheta_q(\beta_1 /\beta_2)\vartheta_q(q^{\lambda-1}\alpha_2/\beta_1)}, \\
&J_{5,1}(x) = -\beta_1^{\lambda-\mu'}x^{\lambda}
\frac{\vartheta_q(\beta_1 x)\vartheta_q(q^{-\lambda}\beta_2 x)\vartheta_q(\beta_2/\alpha_1)}{\vartheta_q(\alpha_2/\beta_1)\vartheta_q(q^{-\lambda}\beta_1\beta_2 x/\alpha_1)\vartheta_q(q^{-\lambda}\alpha_1 x)},\\
&J_{5,3}(x)= - \left( \frac{q}{\beta_1} \right)^{-\lambda+\mu'} x^{\lambda-\mu'}
\frac{\vartheta_q(\beta_1 x)\vartheta_q(\alpha_2 x/q)\vartheta_q(q^{\lambda}\alpha_1/\beta_2)}{\vartheta_q(\beta_1 x/q)\vartheta_q(\alpha_2/\beta_1)\vartheta_q(q^{-\lambda}\beta_1\beta_2 x/\alpha_1)}, \\
&J_{6,1}(x) = -\beta_1^{\lambda-\mu'}x^{\lambda}
\frac{\vartheta_q(q^{\lambda})\vartheta_q(\alpha_1 x)\vartheta_q(\beta_2/\beta_1)}{\vartheta_q(\alpha_2/\beta_1)\vartheta_q(q^{\lambda}\alpha_1/\beta_1)\vartheta_q(q^{-\lambda}\alpha_1 x)},\\
&J_{6,3}(x)= -\left( \frac{\beta_1}{\beta_2} \right)^{\lambda-\mu'}
\frac{\vartheta_q(\beta_2/\beta_1)\vartheta_q(\alpha_2/\beta_2)\vartheta_q(q^{\lambda}\alpha_1/\beta_2)}{\vartheta_q(\beta_1 /\beta_2)\vartheta_q(\alpha_2/\beta_1)\vartheta_q(q^{\lambda}\alpha_1/\beta_1)}.
\end{align*}
The coefficients $K_{m',n'}(x)$\ $(m'=1,3,5,\ n'=1,3)$ are those that replaced $\beta_1$ and $\beta_2$ in $J_{m',n'}(x)$.

\end{document}